\documentclass[a4paper,11pt]{article}
\usepackage{authblk}
\usepackage[english]{babel} 
\usepackage[centertags,intlimits]{amsmath}
\usepackage{amsfonts}
\usepackage{amsthm}
\usepackage{amssymb}
\usepackage{pgf, tikz, pgfplots, float, graphicx, tikz-cd}
\usetikzlibrary{decorations.pathreplacing,calligraphy, arrows}
\usepackage{comment}
\usepackage{cite}

\newtheorem{theorem}{Theorem}[section]
\newtheorem{proposition}[theorem]%
{Proposition}
\newtheorem{definition}{Definition}[section]
\newtheorem{lemma}[theorem]%
{Lemma}
\newtheorem{corollary}[theorem]%
{Corollary}

\numberwithin{equation}{theorem}

\begin{document}

\title{Size of bipartite graphs with given diameter and connectivity constraints}

\author[1,2]{Sonwabile Mafunda}
\affil[1]{Soka University of America, USA}
\affil[2]{University of Johannesburg, South Africa}

\maketitle

\begin{abstract}
In the first part of this paper we determine the maximum size of a (finite, simple, connected) bipartite graph of given order, diameter $d$, and connectivity $\kappa$.\\

It was shown by Ali, Mazorodze, Mukwembi and Vetr\'ik [On size, order, diameter and edge-connectivity of graphs. Acta Math. Hungar. {\bf 152}, (2017)] that for a connected triangle-free graph of order $n$, diameter $d$ and edge-connectivity $\lambda$, the size is bounded from above by about $\frac{1}{4}\left(n-\frac{(\lambda +c) d}{2}\right)^2+O(n)$, where $c\in\{0, \frac{1}{3}, 1\}$ for different values of $\lambda$.

In the second part of this paper we show that this bound by Ali et al. on the size can be improved significantly for a much larger subclass of triangle-free graphs, namely, bipartite graphs of order $n$, diameter $d$ and edge-connectivity $\lambda$. We prove our result only for $\lambda = 2, 3, 4$ because it can be observed from this paper by Ali et al. that for $\lambda\geq 5$, there exists $\ell$-edge-connected bipartite graphs of given order and diameter whose size differs from the maximal size for given minimum degree $\ell$ only by at most a constant. Also, unlike the approach in the proof on the size of triangle-free graphs by Ali et al., our proof employs a completely different technique, which enables us to identify the extremal graphs; hence the bounds presented here are sharp.
\end{abstract}

MSC-class: 05C12 

\section{Introduction}
In this paper we consider only connected, finite bipartite graphs without loops or parallel edges.\\
Let $G$ be a connected, finite graph of order $n \geq 2$ with vertex set $V(G)$ and edge set $E(G)$. The \textit{size} of $G$, denoted $m(G)$ is the cardinality of the edge set of $G$. The \textit{diameter} ${\rm diam}(G)$ is the maximum distance among all pairs of vertices in $G$, i.e. ${\rm diam}(G)=\max\limits_{u,v\in V(G)}d_G(u,v)$, where $d_G(v,u)$ denotes the usual shortest path distance.\\ If $|V(G)|\geq 3$, then the {\it connectivity} $\kappa(G)$  ({\it edge-connectivity} $\lambda(G)$) of $G$ is the smallest number of vertices (edges) whose removal renders $G$ disconnected. The {\it minimum degree} $\delta(G)$ of $G$ is the smallest number of edges incident to a vertex of $G$ among all vertices of $G$.\\

Bounds on the size of a graph in terms of order and diameter were first investigated by Ore \cite{or} in 1968. 
\begin{theorem}{\rm(Ore \cite{or})}\label{thmor}\\
Let $n,d\in\mathbb{N}$. If $G$ is a connected graph of order $n$ and diameter $d$, then
\[m(G)\leq d+\frac{1}{2}\left(n-d-1\right)\left(n-d+4\right).\]
\end{theorem}
In the same paper Ore strengthened this bound by also prescribing connectivity. 
\begin{theorem}{\rm(Ore \cite{or})}\label{thmore}\\
Let $n,d,\kappa\in\mathbb{N}$, with $d \geq 4$ and $\kappa \geq 2$. If $G$ is a $\kappa$-connected graph of order $n$, diameter $d$, then
\[m(G)\leq \frac{1}{2}d\cdot \kappa \left(3\kappa -1\right)-5\kappa^{2}+3\kappa +\frac{1}{2}\left(n-2-\left(d-2\right)\kappa \right)\left(n-3-\left(d-6\right)\kappa\right).\]
\end{theorem}
Later Caccetta and Smyth \cite{ca} and Mukwembi \cite{mu} independently proved a corresponding bounds when connectivity is replaced by the minimum degree $\delta$. Several authors tackled the problem of prescribing edge-connectivity $\lambda$ instead of connectivity or minimum degree. Caccetta and Smyth \cite{ca2} proved sharp upper bounds for $\lambda\geq 8$, leaving a gap for $2\leq\lambda\leq 7$. Asymptotically sharp bounds for these missing values of $\lambda$ were given by Ali et al. \cite{al} and later improved by Dankelmann \cite{da} who proved sharp upper bounds for those values of $\lambda$, and thus completing the determination of the maximum size  of graphs given order, diameter and edge-connectivity. Since the corresponding results for $2\leq \lambda\leq 7$ by Dankelmann \cite{da} and for $\lambda\geq 8$ by Caccetta and Smyth \cite{ca2} all involve cases and are lengthy, we omit them here and refer the reader to the respective papers for brevity. \\
In \cite{al} Ali, Mazorodze, Mukwembi and Vetr\'ik strengthened the bounds on size given either minimum degree $\delta$ or edge-connectivity $\lambda$ by considering graphs containing no triangles. They showed that the bounds they found differ from the actual value by at most $cn$, where $c$ is a constant depending on either $\delta$ or $\lambda$, but not on $d$ or $n$. 
\begin{theorem}{\rm (Ali, Mazorodze, Mukwembi, Vetr\'ik \cite{al})}\label{thm5.7}\\
Let $n,d,\delta\in\mathbb{N}$, with $\delta \geq 2$. If $G$ is a connected triangle-free graph of order $n$, diameter $d$ and minimum degree $\delta$ then
\[m(G)\leq\frac{1}{4}\left[n-\frac{\delta d}{2}\right]^2+O(n).\]
This bound is sharp apart apart from linear term in $n$ if $\delta$ is constant.
\end{theorem}

\begin{theorem}{\rm (Ali, Mazorodze, Mukwembi, Vetr\'ik \cite{al})}\label{thm1.2.14}\\
Let $n,d,\lambda\in\mathbb{N}$, with $\lambda \geq 2$. Let $G$ be a $\lambda$-edge-connected triangle-free graph of order $n$, diameter $d$.
 \[m(G) \leq  \left\{ \begin{array}{c}
\quad\frac{1}{4}\left(n-\frac{(\lambda +1)d}{2}\right)^2+O(n)\quad \textrm{if $\lambda=2, 3$,} \\  
\\
\frac{1}{4}\left(n-\frac{(\lambda +\frac{1}{3})d}{2}\right)^2+O(n)\quad \textrm{if $\lambda=5$,} \\ 
\\
\qquad \frac{1}{4}\left(n-\frac{\lambda d}{2}\right)^2+O(n)\quad \textrm{if $\lambda=4$ or $\lambda\geq 6$.} \\
\end{array} \right. \] 

The bound is sharp apart from a linear term term on the $n$.
\end{theorem}
However, the exact maximum size for a bipartite graph in terms of order~$n$, diameter~$d$ and either connectivity~$\kappa$, edge-connectivity~$\lambda$ or even minimum degree~$\delta$ remained unknown.\\
\\
The goal of this paper is to show that the bounds in Theorem \ref{thmore} and Theorem \ref{thm1.2.14} can be strengthened significantly for bipartite graphs. To achieve this we employ a technique used by Dankelmann in proving similar results on the size of graphs \cite{da}. Our strategy enables us to prove sharp results.

The literature contains several results on the size of graphs ranging from bounds on the size of different classes of graphs to bounds relating the size to other graph parameters. For results on maximum size of a graph given other parameters, for example radius, see \cite{dan, vi}, for domination number see \cite{dank} and for remoteness see \cite{danke}. For results related minimum size given other graph parameters see, for example maximum degree and diameter \cite{er, erd} and minimum degree and diameter \cite{bo}. For graphs of diameter 2, F{\"u}redi~\cite{fu} investigated the maximum size of minimal graphs.

This paper is organised as follows. Section \ref{term&not} introduces the terminology and notation. In Section \ref{section:preliminaries-on-bipartite}, we present and prove results on bipartite graphs that are used in the subsequent sections. Section \ref{vertex} determines the maximum size of a $\kappa$-connected bipartite graph of given order and diameter. Finally, Section \ref{edge} considers edge-connectivity in place of (vertex-)connectivity and provides an upper bound on the size of a $\lambda$-edge-connected graph of given order and diameter.

\section{Terminology and notation}\label{term&not}
We use the following notation.\\
Let $G=(V(G),E(G))$ be a finite, simple, connected graph, where $V(G)$ is a finite nonempty {\em vertex set} and $E(G)$ is a possibly empty {\em edge set}. The cardinalities of the vertex set $|V(G)|$ and edge set $|E(G)|$ are called the {\em order} and {\em size} of $G$ and are denoted by $n(G)$ and $m(G)$ respectively. 

If $v$ is a vertex in $G$, then the {\it neighbourhood} of $v$, denoted by $N(v)$, is the set of all vertices adjacent to $v$. The {\it closed neighbourhood} $N[v]$ of $v$ is the set $N(v)\cup \{v\}$.   

If $G$ is not trivial and not complete, then the {\it connectivity} $\kappa(G)$  ({\it edge-connectivity} $\lambda(G)$) of $G$ is the smallest number of vertices (edges) whose removal renders $G$ disconnected. The {\it minimum degree} $\delta(G)$ of $G$ is the smallest number of edges incident to a vertex of $G$ among all vertices of $G$.

The {\it eccentricity} of $v$, denoted by ${\rm ecc}(v)$, is the distance from $v$ to a vertex farthest from $v$. The {\it diameter} ${\rm diam}(G)$ of $G$ is the largest of all eccentricities of the vertices of $G$. A vertex whose eccentricity equals ${\rm diam}(G)$ is called a {\em peripheral vertex} of $G$. 

The {\it complement} $\overline{G}$ of a graph $G$ is a graph with $V(\overline{G})=V(G)$ and $E(\overline{G})=\{uv\mid u\neq v\in V(G), uv\not\in E(G)\}.$ A graph $G$ is {\it complete} if its complement has no edges. By $K_n$, $\overline{K}_n$, $C_n$ and $P_n$ we mean the complete graph, the empty graph, the cycle, and the path on $n$ vertices. By a {\it triangle} we mean the graph $K_3$. If $F$ is a graph, then we say that $G$ is $F$-{\it free} if $G$ does not contain $F$ as a (not necessarily induced) subgraph. 

We define the {\it sequential sum} of graphs $G=G_1+G_2+\cdots + G_n$ to be the sequential join such that the vertex set 
$V(G)=V(G_1)\cup V(G_2)\cup\cdots\cup V(G_n)$
and the edge set 
$E(G)=E(G_1)\cup E(G_2)\cup\cdots\cup E(G_n)\cup\{uv|\;u\in V(G_i), v\in V(G_{i+1})\quad{\rm where}\quad i=1,2,\ldots, n-1\}$.

If $i\in\mathbb{Z}$, then $N_i(v)$ is the set of all vertices at distance $i$ from $v$, and $n_i$ its cardinality. By $N_{\leq i}(v)$ and $N_{\geq i}(v)$ we mean the set of vertices at distance at most $i$ and at least $i$, respectively, from $v$. Clearly we have $n_i \geq 1$ if and only if $0\leq i\leq {\rm ecc}(v)$. By the {\em distance degree} $X(v)$ of $v$ we mean the sequence $(n_0,n_1,\ldots,n_d)$, where $d$ is the eccentricity of $v$. Here, and for any finite sequence $(n_0,\ldots,n_d)$, we use the convention that $n_i=0$ for all $i\in \mathbb{Z}$ with $i<0$ or $i>d$. With this convention, we define for a finite sequence  $X=(n_0,n_1,\ldots,n_d)$,
\[f(X)=\sum^{d-1}_{i=0}n_in_{i+1}\qquad\mbox{and}\qquad g(X)=\sum^d_{i=0}in_i.\]
If $X$ and $X^{*}$ are two sequences of nonnegative integers, then we say that 
$X^{*}$ {\it beats} $X$ if either $f(X^*)>f(X)$ or $f(X^*)=f(X)$ and $g(X^*)>g(X)$.

We often modify a given sequence $X=(n_0, n_1,\ldots,n_d)$. 
If $a\in \mathbb{N}\cup \{0\}$ and $X^*$ is the sequence obtained from $X$ by adding
$a$ to $n_i$, i.e. if 
$X^*=(n_0,\ldots,n_{i-1}, n_i+a, n_{i+1}\ldots,n_d)$, then we say that $X^*$ is
obtained from $X$ by applying $n_i \leftarrow +a$. Similarly we define the
sequence $X^*$ obtained from $X$ by applying $n_i \leftarrow -a$.
We usually write $X^*$ as $(n_0', n_1',\ldots,n_{d}')$. 

\section{Preliminary results on bipartite graphs}
\label{section:preliminaries-on-bipartite}

In this section we present some results which will be used in the sections that follow.

For a finite sequence $X=(n_0,n_1,\ldots,n_d)$ of positive integers we define
the bipartite graph $G(X)$ by 
\[ G(X) = \overline{K}_{n_0} + \overline{K}_{n_1} +  \ldots + \overline{K}_{n_d}.\]

\begin{definition}\label{def3.1}
{\rm
A finite bipartite graph $G$ of order $n$ and diameter $d$ is said to be \textit{diameter critical} if for any $e\in E(\overline{G})$ we have that $G+e$ is bipartite and $diam(G+e)<diam(G)$.
}
\end{definition}

Since the result below is an analogue of a result due to Ore \cite{or}, for brevity we omit the proof here.

\begin{proposition}\label{prop3.1}
Let $G$ be a connected bipartite graph of order $n$ and diameter $d$. 
Then the following two statements are equivalent. \\
(i) $G+e$ is bipartite and ${\rm diam}(G+e) < {\rm diam}(G)$ for all
    $e\in E(\overline{G})$. \\
(ii) There exist $n_0,n_1,\ldots,n_d \in \mathbb{N}$ 
with $n_0=n_d=1$ and $\sum_{i=0}^d n_i =n$ and
$G = G(X)$.   
\end{proposition} 

In the next two propositions, let $u$ be a vertex of eccentricity $d$ and let $N_i := N_i(u)$ and $n_i := n_i(u)$ for $i \in \mathbb{Z}$.
\begin{proposition}\label{prop3.2}
Let $G$ be a bipartite graph of diameter $d$, let $u$ be a peripheral vertex of $G$, and $X=(n_0,n_1,\ldots,n_{d-1},n_d)$ the distance degree of $u$. Then, for any $x\in N_i$ we have that $N(x)\subseteq N_{i-1}\cup N_{i+1}$.
\end{proposition}

{\bf Proof:} 
It can easily be shown that the neighbourhood of any vertex in $N_i$ lies in $N_{i-1}\cup N_i\cup N_{i+1}$. We show here that in fact if the graph is bipartite then the  neighbourhood of a vertex in $N_i$ lies in $N_{i-1}\cup N_{i+1}$.  It suffices to show that no pair of vertices in $N_i$ are adjacent. Indeed, if there was an edge $xy$ with $x,y\in N_{i}$, then a $u-x$ path of length $i$ together with a $u-y$ path of length $i$ would form a circuit of length $2i+1$, so $G$ has an odd circuit and thus an odd cycle, a contradiction to $G$ being bipartite.\hfill$\Box$\\

Proposition \ref{prop3.3} below shows how the value $f(X)$ of the distance degree of a vertex allows us to bound the size of a bipartite graph. 

\begin{proposition}\label{prop3.3}
Let $G$ be a bipartite graph of diameter $d$, let $u$ be a peripheral vertex of $G$, and $X(u)=(n_0,n_1,\ldots,n_{d-1},n_d)$ the distance degree of $u$. Then
$m(G)\leq f(X)$,
with equality when $G = \overline{K}_{n_0} +  
    \overline{K}_{n_1} +  \ldots +
      \overline{K}_{n_d}$.
\end{proposition}

{\bf Proof:} 
Since $G$ is bipartite, by Proposition \ref{prop3.2} we have that the neighbourhood of any vertex in $N_i$ lies in $N_{i-1}\cup N_{i+1}$. Now, let $x\in N_i$, since $N(x)\subseteq N_{i-1}\cup N_{i+1}$ it follows that \[{\rm deg}(x)\leq |N_{i-1}\cup N_{i+1}|=n_{i-1}+n_{i+1}.\] Now consider the neighbourhood of all $x\in N_i$, we have that,
\[\sum_{x\in N_i}{\rm deg}(x)\leq n_i(n_{i-1}+n_{i+1})\] Since $i=0,1,\ldots, d$ we have that for all $i$, the degree sum of $G$ is 
\[\sum_{i=0}^{d}\sum_{x\in N_i}{\rm deg}(x)\leq \sum_{i=0}^{d}n_i(n_{i-1}+n_{i+1}).\] Since $n_{-1}=n_{d+1}=0$, we have that 
\begin{align*}
\sum_{i=0}^{d}\sum_{x\in N_i}{\rm deg}(x)&\leq 2n_0n_1+ 2n_1n_2+\ldots +2n_{d-1}n_d\\ &=2\sum_{i=0}^{d-1}n_in_{i+1}.
\end{align*} 
But $\sum_{i=0}^{d-1}\sum_{x\in N_i}{\rm deg}(x)$ is the degree sum of $G$. And by the Handshaking Lemma we know that
\begin{align*}
m(G)&=\frac{1}{2}\sum_{i=0}^{d-1}\sum_{x\in N_i}{\rm deg}(x)\\
&\leq\sum_{i=0}^{d-1}n_in_{i+1}=f(X)
\end{align*}
as desired. 
For equality, consider $G = \overline{K}_{n_0} + \overline{K}_{n_1} +  \ldots + \overline{K}_{n_d}$, this is a bipartite graph where every $x\in N_i$ is adjacent to every vertex in $N_{i+1}$ and vice versa.\\
Clearly it follows that $\sum_{x\in N_i}{\rm deg}(x)=n_i(n_{i-1}+n_{i+1})$, which implies that $\sum_{i=0}^{d}\sum_{x\in N_i}{\rm deg}(x)=2\sum_{i=0}^{d-1}n_in_{i+1}$. And now, with similar reasoning as above we have that $m(G)=f(X)$. And clearly, if $m(G)=f(X)$, then every vertex in $N_i$ has degree $n_{i-1}+n_{i+1}$, and is thus adjacent to every vertex in $N_{i-1}\cup N_{i+1}$, for every $i\in \{0,1,\ldots,d\}$. Hence we have $G=\overline{K}_{n_0} +  
    \overline{K}_{n_1} +  \ldots +
      \overline{K}_{n_d}$, which completes the proof.
\hfill$\Box$\\

While our goal is to maximise the size, and thus maximise $f$, it turns out that often there are many sequences that maximise the value of $f$. Considering also $g$ allows us to determine a unique sequence that maximises $f$, and among those maximising $f$, maximises $g$. 

\section{Maximum size of $\kappa$-connected bipartite graphs}\label{vertex}
In this section we give a sharp upper bound on the size of a 
bipartite graph of order $n$, diameter $d$ and connectivity $\kappa$. We prove our bound by first demonstrating some properties of the distance degree $X(v)$ of an arbitrary vertex $v$ of a bipartite graph, and then determining a sequence that maximises the function $f$ among all sequences with these properties. 

\begin{proposition} \label{prop3.4}
Let $u$ be a vertex of a connected, bipartite graph of order $n$, diameter $d\geq 3$ and connectivity $\kappa$. Let $X(u)=(n_0,n_1,\ldots,n_d)$. Then the following hold: 
\begin{enumerate}
\item[(A1)] $n_0=1$,
\item[(A2)] $\sum\limits_{i=0}^d n_i=n$,
\item[(A3)] if $i\geq 1$ and $n_i\geq 1$, then $n_1, n_2, \cdots, n_{i-1}\geq 1$,
\item[(A4)] $n_i\geq\kappa$, for all $i=1,2,\ldots, d-1$,
\item[(A5)] $n_d=1$.
\end{enumerate}
\end{proposition}

{\bf Proof:} 
Let $N_i:=N_i(u)$ and $n_i:=n_i(u)$. 
Clearly, (A1), (A2) and (A3) hold. \\
(A4): For $1\leq i\leq d-1$; suppose to the contrary that there exists $1\leq i\leq d-1$ such that $n_i\leq \kappa-1$. Now, removal of all $n_i$ vertices disconnects $G$,
since there is no path from any $v\in N_{\leq i-1}[u]$ to any vertex $x\in N_{\geq i+1}(u)$, thus a contradiction to the claim that at least $\kappa$ vertices must be removed for $G$ to be disconnected.
Therefore, if $G$ is $\kappa$-connected, then $n_{_i}\geq \kappa$ for all $i=1,2,\ldots,d-1$, which is (A4). \\
(A5): Suppose to the contrary $n_d>1$. Let $v$ be a vertex at distance $d-2$ from $u$ and both $x$ and $y$ be vertices of distance $d$ from $u$.
Clearly, $e=vx$ is not an edge in $G$ by Proposition \ref{prop3.2}.  Now, adding $e$ in $G$ we get that $G+e$ is bipartite and ${\rm diam}(G+e)={\rm diam}(G)$, a contradiction to the claim that $G$ is among graphs of maximum size with order $n$ and diameter $d$. Hence, $n_d=1$, which is (A5).
\hfill$\Box$\\

\begin{definition}\label{def3.2}
{\rm
For $n,d, \kappa \in \mathbb{N}$ with $d\geq 3$ and $\kappa\geq 2$, 
define the finite sequence $X^{\kappa}_{n,d}$ to be the sequence of length $d+1$ given by
 \[X_{n,d}^\kappa =  \left\{ \begin{array}{c}
\left(1,\big\lfloor\frac{n-2}{2}\big\rfloor,\big\lceil\frac{n-2}{2}\big\rceil,1\right)\\
\textrm{if $d=3$ and $n\geq 2\kappa +2$,} \\  
\\
\left(1, \kappa, \kappa, \kappa, \ldots, \kappa, \big\lfloor\frac{n-(d-3)\kappa-2}{2}\big\rfloor, \big\lceil\frac{n-(d-3)\kappa-2}{2}\big\rceil, \kappa, 1\right)\\
\textrm{if $d\geq 4$ and $n\geq (d-1)\kappa +2$.} \\  
\end{array} \right. \] 
}
\end{definition}

We now show that, for given $n$, $d$ and $\kappa$, a sequence that satisfies (A1) - (A5), and that is not beaten by any other such sequence, must equal $X_{n,d}^{\kappa}$.\\

\begin{lemma}\label{lem3.6}
Let $n,d, \kappa\in \mathbb{N}$ with $d\geq 4$, $\kappa\geq 2$ and $n\geq (d-1)\kappa +2$. If $X=(n_0,\ldots,n_d)$ is a sequence satisfying {\rm (A1) - (A5)} that is not beaten by any other sequence, then $X=X_{n,d}^{\kappa}$.
\end{lemma}

{\bf Proof:} 
Let $X$ be a sequence satisfying (A1) - (A5). It suffices to prove the lemma for a sequence which is not beaten by any other sequence satisfying (A1) - (A5), so let $X$ be such a sequence. Let $n$, $d$ and $\kappa$ be fixed. By (A1) and (A5) we know that $n_0=n_d=1$. Now in three claims we prove properties of the sequence $X$ by showing that otherwise we can modify $X$ to obtain a sequence 
$X^*= (n_0',n_1',\ldots,n_{d}')$ which satisfies (A1) - (A5) but beats $X$, thus obtaining a contradiction.\\[1mm]
{\sc Claim 1:} If $i=1, 2,\ldots, d-4$ then $n_i=\kappa$. \\
Since, $n_{i}\geq \kappa$ for all $i=1,2,\ldots, d-4$ by (A4) we are only left to show that $n_{i}\leq \kappa$. Suppose to the contrary that for some $i\in \{1, 2,\ldots, d-4\}$, we have $n_i\geq\kappa +1$. Let $j$ be the smallest integer such that $n_j\geq\kappa +1$ and $1 \leq j \leq d-4$. Let $X^*$ be the sequence obtained from $X$
by applying $n_j\leftarrow-1$ and $n_{j+1}\leftarrow+1$,
clearly $X^*$ satisfies all conditions (A1) - (A5). Then $f(X^*)\geq f(X)$ and $g(X^*)=g(X)+1$. It follows that the sequence $X^*$ beats sequence $X$, a contradiction. Hence $n_i\leq \kappa$, which implies $n_i=\kappa$ for all $i\in \{1,2,\ldots, d-4\}$, as desired.\\[1mm]
{\sc Claim 2:} $n_{d-1}=\kappa$. \\
By (A4) we have that $n_{d-1}\geq \kappa$. We claim that $n_{d-1}$ is not greater than $\kappa$. Suppose to the contrary that it is. We know from Claim 1 that $n_{d-4} =\kappa$ which is greater than $1$, hence if $X^*$ is the sequence obtained from $X$ by applying $n_{d-1}\leftarrow-1$ and $n_{d-3}\leftarrow+1$. Then clearly $X^*$ satisfies (A1) - (A5) and since $f(X^*)=f(X)+(\kappa -1)>f(X)$ we have that $X^*$ beats $X$, a contradiction. Therefore $n_{d-1}=\kappa$.\\[1mm]
{\sc Claim 3:} $n_{d-3}=\big\lfloor\frac{n-(d-3)\kappa-2}{2}\big\rfloor$ and $n_{d-2}=\big\lceil\frac{n-(d-3)\kappa-2}{2}\big\rceil$.\\
Since $n_{d-3}, n_{d-2}\geq \kappa$ by (A4) we have that $n_{d-3}+n_{d-2}\geq 2\kappa$. And since $n_{d-3}+n_{d-2}=n-(d-3)\kappa -2$ we have that $n_{d-3}+n_{d-2}=n-(d-3)\kappa -2$ must be at least $2\kappa$. Now we show that $0\leq n_{d-2}-n_{d-3}\leq 1$.\\
Clearly, $n_{d-2}-n_{d-3}\geq 0$ since otherwise, if $n_{d-3}>n_{d-2}$, swapping the values of $n_{d-3}$ and $n_{d-2}$ would yield a sequence that satisfies (A1) - (A5) and has the same $f$-value, but a greater $g$-value, a contradiction. Hence it is only left to show that $n_{d-2}-n_{d-3}\leq 1$.\\
Suppose to the contrary that $n_{d-2}-n_{d-3}\geq 2$. Let $X^*$ be the sequence obtained from $X$ by applying $n_{d-2}\leftarrow-1$ and $n_{d-3}\leftarrow+1$. Then clearly $X^*$ satisfies (A1) - (A5) and since  $f(X^*)=f(X)+n_{d-2}-(n_{d-3}+1)>f(X)$, we have that $X^*$ beats $X$, a contradiction. Hence $n_{d-2}-n_{d-3}\leq 1$, which proves the claim, thus completes the proof of the lemma.
\hfill$\Box$\\

\noindent Recall that for a finite sequence $X=\left(x_0, x_1,\cdots, x_d\right)$ of positive integers we define the bipartite graph $G(X)$ by \[G(X)=\overline{K}_{x_0} +\overline{K}_{x_1}+\cdots +\overline{K}_{x_d}.\]

\begin{theorem}\label{thm3.7}
Let $G$ be a $\kappa$-connected, bipartite graph of order $n$ and diameter $d$, where $d\geq 4$, $\kappa \geq 2$ and $n\geq (d-1)\kappa +2$. Then 
\[ m(G) \leq m(G(X_{n,d}^{\kappa})). \]
\end{theorem}

{\bf Proof:} 
Let $n, d$ and $\kappa$ be fixed. Assume $G$ is a $\kappa$-connected bipartite graph of order $n$ and diameter $d$, and $u$ is a vertex of eccentricity $d$. Let $X(u)$ be the distance degree of $u$. It follows by Proposition \ref{prop3.3} that 
\[m(G) \leq f(X(u)).\]
Since $X(u)$ satisfies (A1) - (A5), and since by Lemma \ref{lem3.6} sequence $X(u)$ does not beat $X_{n,d}^{\kappa}$, we have 
\[ f(X) \leq f(X_{n,d}^{\kappa})\quad{\rm for}\quad d\geq 4.\]
Similarly, from Lemma \ref{lem3.5}, we have
\[ f(X) \leq f(X_{n,3}^{\kappa}).\]
If $u$ is the vertex of $G(X_{n,d}^{\kappa})$ contained in $K_{x_0}$, then 
the distance degree of $u$ is $X_{n,d}^{\kappa}$, so 
\[ f(X_{n,d}^{\kappa}) \leq m(G(X_{n,d}^{\kappa})). \] 
Combining the all these (in)equalities yields the theorem.
\hfill$\Box$\\

\begin{theorem}\label{lem3.5}
Let $n,d, \kappa\in \mathbb{N}$ with $d=3$, $\kappa\geq 2$ and $n\geq 2\kappa +2$. Let $G$ be a $\kappa$-connected bipartite graph of order $n$ and diameter $d$ that has maximum size among all such graphs. Then, 
 \[m(G) \leq m\left(\overline{K}_1 + \overline{K}_{\big\lfloor\frac{n-2}{2}\big\rfloor} + \overline{K}_{\big\lceil\frac{n-2}{2}\big\rceil} + \overline{K}_1\right).\]
\end{theorem}

{\bf Proof:} 
Let $v$ be a peripheral vertex of $G$ and $X(v)=(n_0,n_1,n_2,n_3)$ a distance degree of $v$. Then we have that $G(X)$ is bipartite. We assume by Proposition \ref{prop3.3} that $G=G(X)$. Hence, it should suffice to determine the integers $n_0, n_1, n_2, n_3$. Following from Proposition \ref{prop3.4} we have that $\overline{K}_{n_0}=\overline{K}_1$ and $\overline{K}_{n_3}=\overline{K}_1$. Since by (A4) we have $n_1, n_2\geq \kappa$, we have $n \geq 2\kappa +2$. We also have $0 \leq n_2-n_1 \leq 1$. Indeed, if $n_2-n_1 <0$ then it is easy to verify that the size of $\overline{K}_1 + \overline{K}_{n_1-1} + \overline{K}_{n_2+1} + \overline{K}_1$ is not less than the size of $\overline{K}_1 + \overline{K}_{n_1} + \overline{K}_{n_2} + \overline{K}_1$. Similarly, if $n_2-n_1 >1$ then it is easy to verify that the size of $\overline{K}_1 + \overline{K}_{n_1+1} + \overline{K}_{n_2-1} + \overline{K}_1$ is not less than the size of $\overline{K}_1 + \overline{K}_{n_1} + \overline{K}_{n_2} + \overline{K}_1$, which completes the proof of the theorem.
\hfill$\Box$\\

\noindent Since $G(X_{n,d}^{\kappa})$ is a bipartite graph of order $n$, diameter $d$ and connectivity $\kappa$, the bounds obtained in Theorem \ref{thm3.7} and Theorem \ref{lem3.5} are sharp. Evaluating the size of $G(X_{n,d}^{\kappa})$ yields the following corollary.\\

\begin{corollary}\label{cor3.8}
Let $n, d, \kappa\in \mathbb{N}$, with $\kappa \geq 2$ and $d\geq 3$. If $G$ is a $\kappa$-connected bipartite graph of order $n$ and diameter $d$, then
 \[ m(G)\leq  \left\{ \begin{array}{c}
\frac{n^2-4}{4},\\
\textrm{if $d=3$, $n\geq 2\kappa +2$ and $n$ is even,} \\  
\\
\frac{n^2-5}{4},\\
\textrm{if $d=3$, $n\geq 2\kappa +2$ and $n$ is odd,} \\  
\\
\left(n-2\kappa\right)\kappa +\big\lceil\frac{n-(d-3)\kappa-2}{2}\big\rceil \big\lfloor\frac{n-(d-3)\kappa-2}{2}\big\rfloor,\\
\textrm{if $d\geq 4$ and $n\geq (d-1)\kappa +2$} \\  
\end{array} \right. \] 
and this bound is sharp.
\end{corollary}

\section{Maximum size of $\lambda$-edge-connected bipartite graphs}\label{edge}
In this section we give an upper bound on the size of a bipartite graph in terms of order, diameter and edge-connectivity, which improves on the bound by Ali et al. \cite{al} on maximum size of a triangle-free graph in terms of order, diameter and edge-connectivity, the bounds proved here are sharp.

Our strategy of proving the bound on maximum size in terms of order, diameter and edge-connectivity follows a similar method as that in Section \ref{vertex}.

\begin{proposition} \label{prop4.4}
Let $u$ be a vertex of a connected, bipartite graph of order $n$, diameter $d\geq 3$ and edge-connectivity $\lambda$. Let $X(u)=(n_0,n_1,\ldots,n_d)$. Then the following hold: 
\begin{enumerate}
\item[(B1)] $n_0=1$,
\item[(B2)] $\sum\limits_{i=0}^d n_i=n$,
\item[(B3)] if $i\geq 1$ and $n_i\geq 1$, then $n_1, n_2, \cdots, n_{i-1}\geq 1$,
\item[(B4)] $n_in_{i+1}\geq\lambda$, for all $i=0,1,\ldots, d-1$,
\item[(B5)] $n_{i-1}+n_{i+1}\geq\lambda$, for all $i=0,1,\ldots, d$,
\item[(B6)] $n_d=1$.
\end{enumerate}
\end{proposition}

{\bf Proof:} 
Let $N_i:=N_i(u)$ and $n_i:=n_i(u)$. 
(B1), (B2), (B3) and (B6) are the properties (A1), (A2), (A3) and (A5), respectively, in 
Proposition \ref{prop3.4}. \\
(B4): Let $E_{i,i+1}$ be the set of all edges joining $N_{\leq i}[u]$ to $N_{\geq i+1}(u)$.
Then every edge in $E_{i,i+1}$ joins a vertex in $N_i$ to a vertex in $N_{i+1}$, so $|E_{i,i+1}|\leq n_in_{i+1}$. Clearly $E_{i,i+1}$ is an edge-cut in $G$, which implies $|E_{i,i+1}| \geq \lambda$. Hence, we have that $\lambda\leq |E_{i,i+1}| \leq n_in_{i+1}$. This implies $\lambda\leq n_in_{i+1}$, which proves (B4).\\
(B5): Since $G$ is bipartite then by Proposition \ref{prop3.2} $N(x)\subseteq N_{i-1}\cup N_{i+1}$, for any $x\in N_i$.
This implies that a vertex $x\in N_i$ has degree, ${\rm deg}(x)\leq n_{i-1}+n_{i+1}$. But $G$ is $\lambda$-edge connected, hence $n_{i-1}+n_{i+1}\geq\lambda$.
\hfill$\Box$\\

\begin{definition}\label{def4.1}
{\rm
For $n,d, \lambda \in \mathbb{N}$ with $\lambda =2,3,4$, $d\geq 6$ and 
define the finite sequence $X^{\lambda}_{n,d}$ to be the sequence of length $d+1$ given 
\[ X^{\lambda}_{n,d} = \left\{ \begin{array}{c}
   \left(1, 2, 1,2,1,2,\dots,1,2, \Big\lfloor \frac{2n-3(d-1)}{4} \Big\rfloor,  \Big\lceil \frac{2n-3(d-1)}{4} \Big\rceil,2,1\right)  \\
         \textrm{if $\lambda = 2$, $d$ odd and $n \geq \frac{3}{2}(d+1)$,} \\  
         \\
                 \left(1, 2, 1,2,1,2,\dots,1,2, 1, 2,1\right)\quad\mbox{where}\quad n-\frac{3}{2}d+2=3, \\
        \left(1, 2, 1,2,1,2,\dots,1,2, \Big\lceil\frac{2n-3d-4}{4}\Big\rceil, \Big\lfloor\frac{2n-3d-4}{4}\Big\rfloor,1\right)\quad\mbox{otherwise},  \\
        \textrm{if $\lambda = 2$, $d$ even and $n \geq \frac{3}{2}(d+1)$,} \\ 
        \\
   \left(1, \lambda, \lambda-1,2,2,2,\ldots, 2,2,n_{d-3}, n_{d-2},\lambda,1\right)  \\
         \textrm{if $\lambda = 3,4$ and $n\geq 3\lambda+5$,}             
      \end{array} \right. \] 
 where the values for $n_{d-3}$ and $n_{d-2}$ when $d\geq 6$ and $\lambda =3$ are as follows:
 \[ \left(n_{d-3}, n_{d-2}\right) = \left\{ \begin{array}{c}
    \left(2,2\right) \\
         \textrm{if $n-2d-3\lambda+11=4$,} \\  
         \\
   \left(\frac{n-2d+1}{2}, \frac{n-2d+3}{2}\right) \\
         \textrm{if $n-2d-3\lambda+11$ is odd and $n-2d-3\lambda+11\geq 5$,} \\  
         \\
            \left(\frac{n-2d}{2}, \frac{n-2d+4}{2}\right) \\
         \textrm{if $n-2d-3\lambda+11$ is even and $n-2d-3\lambda+11\geq 6$,}          
      \end{array} \right. \] 
      
       and the values for $n_{d-3}$ and $n_{d-2}$ when $d\geq 6$ and $\lambda =4$ are as follows:
 \[ \left(n_{d-3}, n_{d-2}\right) = \left\{ \begin{array}{c}
         \left(2, 3\right) \\
         \textrm{if $n-2d-3\lambda+11=5$,} \\  
         \\
   \left(\frac{n-2d-20}{2}, \frac{n-2d+2}{2}\right)\\
         \textrm{if $n-2d-3\lambda+11$ is odd and $n-2d-3\lambda+11\geq 7$,} \\  
         \\
            \left(\frac{n-2d-3}{2}, \frac{n-2d+1}{2}\right)\\
         \textrm{if $n-2d-3\lambda+11$ is even and $n-2d-3\lambda+11\geq 6$.}           
      \end{array} \right. \]     
}
\end{definition}

We now show that, for given $n$, $d$ and $\lambda$, with $d\geq 6$ and $\lambda=3,4,5$ a sequence that satisfies (B1) - (B6), 
and that is not beaten by any other such sequence, must equal $X_{n,d}^{\lambda}$. 

\begin{lemma}\label{lem4.6}
Let $n,d, \lambda\in \mathbb{N}$ with $d\geq 6$, $\lambda = 2,3,4$ and $n\geq (d-1)\kappa +2$. If $X=(n_0,\ldots,n_d)$ is a sequence satisfying {\rm (B1) - (B6)} that is not beaten by any other sequence, then $X=X_{n,d}^{\lambda}$.
\end{lemma}

{\bf Proof:} 
Let $X$ be a sequence satisfying (B1) - (B6). It suffices to prove the lemma for a sequence which is not beaten by any other sequence satisfying (B1) - (B6), so let $X$ be such a sequence. Let $n$, $d$ and $\lambda$ be fixed. By (B1) and (B6) we know that $n_0=n_d=1$. Now in a sequence of claims we prove properties of the sequence $X$ by showing that otherwise we can modify $X$ to obtain a sequence 
$X^*= (n_0',n_1',\ldots,n_{d}')$ which satisfies (B1) - (B6) but beats $X$, thus obtaining a contradiction.\\[1mm]
{\sc Claim 1:} If $i\leq d-3$ and $n_{i-1}=1$ then $n_i=\lambda$. \\
Since, $n_{i-1}=1$, condition (B4) for $i\leq d-3$ yields $n_{i}\geq \lambda$. So we only need to show is $n_{i}\leq \lambda$. Suppose to the contrary that $n_{i}> \lambda$. Now, if $X^*$ is the sequence obtained from $X$ by applying $n_i\leftarrow-(n_i-\lambda)$ and $n_{i+2}\leftarrow+(n_i-\lambda)$ clearly, $X^{*}$ satisfies all conditions (B1) - (B6).
Then $f(X^*)=f(X)+(n_i-\lambda)(n_{i+3}-1)\geq f(X)$ and $g(X^*)=g(X)+2(n_{i}-\lambda)>g(X)$. It follows that the sequence $X^{*}$ beats sequence $X$, a contradiction. Hence $n_i\leq\lambda$, which implies, $n_i=\lambda$, as desired.\\[1mm]
{\sc Claim 2:} If $\lambda=2$, $i\leq d-4$ and $n_{i-1}=\lambda$ then $n_i=1$.\\
By (B4) we have that $n_i\geq 1$, hence we only left to show that $n_i\leq 1$. Suppose to the contrary that for some integer $i\in \{0,1,\ldots, d-4\}$ we have $n_{i-1}=\lambda$ and $n_{i}\geq 2$. Let $j$ be one such $i$.\\
{\sc Case 2.1:} $n_{j+1}=1$.\\
Then $n_{j+2}=\lambda$ by Claim 1. Now, let $X^*$ be the sequence obtained from $X$ by applying $n_j\leftarrow-(n_j-1)$ and $n_{j+1}\leftarrow+(n_j-1)$. Then clearly $X^*$ satisfies (B1), (B2), (B3) and (B6). We now show that $X^*$ satisfies (B4). We have that (B4) holds for all $i$ except possibly when $i=j-1, j, j+1$ since only the values $n_j$ and $n_{j+1}$ have changed. Since $i\leq d-4$ and $n_{j-1}'=n_{j-1}=\lambda$, $n_{j}'=n_{j}-(n_j-1)=1$, $n_{j+1}'=n_{j+1}+n_j-1\geq 2$, $n_{j+2}'=n_{j+2}=\lambda$, it is clear that (B4) holds for $i=j-1, j, j+1$, hence holds for all $i$. Next, we show that (B5) holds.  Clearly (B5) holds for all $i$ except possibly when $i=j-1, j, j+1, j+2$. Since $n_0=1$ and $n_1=\lambda$ by Claim 1 we have that $2\leq j\leq d-4$. Now, by (B4) and since $\lambda =2$, it is clear that (B5) holds for all $i$. Since $f(X^*)=f(X)+(n_{j}-1)(n_{j+2}+1-n_{j-1}-n_{j+1})\geq f(X)$ and $g(X^*)=g(X)+(n_j-1)>g(X)$. This implies that sequence $X^{*}$ beats sequence $X$, a contradiction. Hence, if $i\leq d-4$ and $n_{i-1}=\lambda$ then $n_i=1$.\\
{\sc Case 2.2:} $n_{j+1}\geq 2$.\\
{\sc Case 2.2.1} $n_{j+3}=1$.\\
Then $n_{j+4}=\lambda$ by Claim 1 and $n_{j+2}\geq 2$ by (B4). It is clear that $n_{j+2}\leq 2$. This follows by the assumption that if $n_{j+2}$ was greater than $2$ then the operation $n_{j+2}\leftarrow -(n_{j+2}-2)$ and $n_{j}\leftarrow +(n_{j+2}-2)$ would produce a sequence that beats $X$. Hence, $(n_{j+2}, n_{j+3}, n_{j+4})=(2,1,2)$. Now, let $X^*$ be the sequence obtained from $X$ by simultaneously applying
$n_j\leftarrow-\left(n_j-1\right)$ and $n_{j+1}\leftarrow -\left(n_{j+1}-2\right)$, 
$n_{j+2}\leftarrow +\left(n_{j+1}-2\right)$ and $n_{j+3}\leftarrow +\left(n_j-1\right)$. Then $X^*$ satisfies (B1), (B2) and (B3). We now show that $X^*$ satisfies (B4). We have that (B4) holds for all $i$ except possibly when $i=j-1, j, j+1,j+2,j+3$ since only the values $n_j$, $n_{j+1}$, $n_{j+2}$ and $n_{j+3}$ have changed. Since $2\leq i\leq d-4$ and $n_{j-1}'=n_{j-1}=\lambda$, $n_{j}'=n_{j}-(n_j-1)=1$ and $n_{j+1}'=n_{j+1}-(n_{j+1}- 2)=2$ it  follows that (B4) holds for $i=j-1, j$. Now, for $i=j+1, j+2, j+3$ we note that $n_{j+2}'=n_{j+2}+\left(n_{j+1}-2\right) =n_{j+1}\geq \lambda$, $n_{j+3}'=n_{j+3}+ \left(n_{j}-1\right)=n_{j}\geq  \lambda$, it is clear that (B4)  holds for $i=j+1, j+2, j+3$, hence holds for all $i$. Next, we show that (B5) holds.  Clearly (B5) holds for all $i$ except possibly when $i=j-1, j, j+1, j+2, j+3, j+4$. Since $n_0=1$ and $n_1=\lambda$ by Claim 1 we have that $2\leq j\leq d-4$. Now, by (B4) and since $\lambda =2$, it is clear that (B5) holds for all $i$ except possibly for $i=j+4$. But we know that $n_{j+3}'=n_{j+3}+\left(n_{j}-1\right)= n_{j}\geq \lambda$, hence (B5) holds for $i=j+4$, thus (B5) holds for all $i$. Since $f(X^*)=f(X)$ and $g(X^*)=g(X)+3n_j+n_{j+1}-5>g(X)$ it then follows that the sequence $X^{*}$ beats sequence $X$, a contradiction. Hence, if $i\leq d-4$ and $n_{i-1}=\lambda$ then $n_i=1$.\\
{\sc Case 2.2.2} $n_{j+3}\geq 2$.\\
Let $X^*$ be the sequence obtained from $X$ by applying $n_j\leftarrow-(n_j-1)$ and $n_{j+2}\leftarrow+(n_j-1)$. Then clearly $X^*$ satisfies (B1), (B2), (B3) and (B5). We now show that $X^*$ satisfies (B4). We have that (B4) holds for all $i$ except possibly when $i=j-1, j, j+1, j+2$ since only the values $n_j$ and $n_{j+2}$ have changed. Since $i\leq d-4$ and $n_{j-1}'=n_{j-1}=\lambda$, $n_{j}'=n_{j}-(n_j-1)=1$, $n_{j+1}'=n_{j+1}\geq 2$, $n_{j+2}'=n_{j+2}+n_j-1\geq 2$, it is clear that (B4)  holds for $i=j-1, j, j+1$, hence holds for all $i$. Next, we show that (B5) holds.  Clearly (B5) holds for all $i$ except possibly when $i=j-1, j+1, j+3$. Since $n_0=1$ and $n_1=\lambda$ by Claim 1 we have that $2\leq j\leq d-4$. Now, by (B3) and since $\lambda =2$, it is clear that (B5) holds for all $i$. Since, $f(X^*)=f(X)+(n_{j+3}-2)(n_{j}-1)\geq f(X)$ and $g(X^*)=g(X)+2(n_j-1)>g(X)$, it then follows that the sequence $X^{*}$ beats sequence $X$, a contradiction. Hence, if $i\leq d-4$ and $n_{i-1}=\lambda$ then $n_i=1$, as desired.\\[1mm]
{\sc Claim 3:} If $\lambda=3,4$, $n_i=2$ and $n_{i+1}\geq 3$, then $n_{i+2}\geq n_{i+1}\geq 3$, $n_{i+3}=\lambda$, $n_{i+4}=1$ and $i=d-4$.\\
Assume $n_{i}, n_{i+1},\dots , n_{i+4}$ and let $n_{i}=2$, $n_{i+1}\geq 3$. Then $n_{i+2}\geq 3$, otherwise the operations $n_{i+1}\leftarrow-1$ and $n_{i+2}\leftarrow+1$ on $X$ produces a sequence that satisfies conditions (B1) - (B6) and beats sequence $X$. Similarly, $n_{i+4}=1$, otherwise the operations $n_{i+1}\leftarrow-1$ and $n_{i+3}\leftarrow+1$ on $X$ produces a sequence that clearly satisfies conditions (B1) - (B4) and (B6). Also, since $n_{i}=2$, it follows that $n'_{i-1}\geq 2$, which implies (B5) holds for $i$ and, for $i+1$ and $i+4$ (B5) holds since $n'_{i+3}>n_{i+3}$. Hence  (B5) holds for all integers $i\in\{1,2,\ldots, d\}$ in the sequence resulting from the above operations and this sequence would beat sequence $X$. Now, since $n_{i}=2$, $n_{i+1}, n_{i+2}\geq 3$ and $n_{i+4}=1$, it is clear by condition (B5) that $n_{i+3}\geq \lambda$. So, we only need to show that $n_{i+3}\leq \lambda$. Now, suppose to the contrary that $n_{i+3}\geq \lambda +1$. Let $X^*$ be the sequence obtained from $X$ by applying $n_{i+3}\leftarrow-(n_{i+3}-\lambda)$ and $n_{i+1}\leftarrow+(n_{i+3}-\lambda)$. Clearly, $X^*$ satisfies (B1), (B2), (B3) and (B6). We now prove that $X^*$ satisfy (B4) and (B5). Clearly (B4) holds for all integers $0,1,\ldots, d-1$ except possibly $i, i+1, i+2, i+3$ since only the values $n_{i+1}$ and $n_{i+3}$ have changed. Since $n'_{i+1}>n_{i+1}$ and $n_i, n_{i+2}$ remained unchanged, we have that (B4) holds for $i$ and $i+1$. Similarly, since $n'_{i+3}=\lambda$ and $n_{i+2}, n_{i+4}\geq 1$, we have that (B4) holds for $i+2$ and $i+3$ as well. We now show that (B5) holds. It is clearly that (B5) holds for all integers $0,1,\ldots, d$ since the only possible integers where (B5) wouldn't hold are $i, i+2, i+4$. But, $n'_{i+1}>n_{i+1}$, which implies that $n'_{i-1}+n'_{i+1}>n_{i-1}+n_{i+1}$, also, $n'_{i+3}=\lambda$, which implies that for any nonnegative integers $n'_{i+1}$ and $n'_{i+5}$ we have that $n'_{i+1}+n'_{i+3}, n'_{i+3}+n_{i+5}\geq \lambda$. Hence, (B5) holds for all integers $0,1,\ldots, d$. And $f(X^*)=f(X)+(n_{i+3}-\lambda)> f(X)$, which implies that the sequence $X^{*}$ beats sequence $X$, a contradiction. Hence, $n_{i+3}\leq \lambda$, thus, $n_{i+3}= \lambda$. Now, it is clear that $n_{i+2}\geq n_{i+1}+\lambda -2$, otherwise the operations $n_{i+1}\leftarrow-1$ and $n_{i+2}\leftarrow+1$ on $X$ produces a sequence that satisfies all conditions (B1) - (B6) and beats sequence $X$.\\ Finally, we show that $i=d-4$. Suppose not, then, since $n_{i+4}=1$ we know by Claim 1 that $n_{i+5}=\lambda$. But simultaneously applying $n_{i+1}\leftarrow-1$ and $n_{i+3}\leftarrow+1$, $n_{i+2}\leftarrow-1$ and $n_{i+4}\leftarrow+1$ on $X$ we have that the resulting sequence satisfies all conditions (B1) - (B6) and beats sequence $X$, hence $i=d-4$. This completes the proof for this claim.\\[1mm]
{\sc Claim 4:} If $\lambda=3,4$ and $n_{i}=1$, then $i=0$ or $i=d$.\\
Suppose to the contrary that there exists some integer $2\leq i\leq d-1$ such that $n_{i}=1$. Let $j$ be the smallest such $i$. We know by Claim 1 that $n_{j+1}=\lambda$ and by condition (B4) that $n_{j-1}\geq \lambda$. Clearly, by conditions (B5) $j+1\neq d$ and $n_{j+2}\geq \lambda-1$. Also, $j+2\neq d$, otherwise the operations $n_{j+2}\leftarrow-(n_{j+2}-1)$ and $n_{j}\leftarrow+(n_{j+2}-1)$ on $X$ produces a sequence that satisfies all conditions (B1) - (B6) and beats sequence $X$. Consider the following cases depending on $j$. By conditions (B1) and (B5) and Claim 1 we have that $(n_0,n_1,n_{2})=(1,\lambda,n_{2})$, where $n_{2}\geq \lambda -1$.\\
{\sc Case 4.1:} $j=3,4$.\\
Since $n_{j-1}\geq \lambda$, we have that $n_{j-2}\leq \lambda$, otherwise the operations $n_{j-2}\leftarrow-1$ and $n_{j}\leftarrow +1$ on $X$ produces a sequence that satisfies all conditions (B1) - (B6) and beats sequence $X$. Now, if $X^*$ is the sequence obtained from $X$ by applying $n_{j-1}\leftarrow-1$ and $n_{j}\leftarrow+1$ clearly, $X^{*}$ satisfies all conditions (B1) - (B6). Then $f(X^*)=f(X)+(n_{j-1}-n_{j-2})+(\lambda -2)> f(X)$, implying that the sequence $X^{*}$ beats sequence $X$, a contradiction.\\
{\sc Case 4.2:} $j\geq 5$.\\
Since $i=j$ is the smallest integer greater than $2$ with $n_i=1$, we have that $n_{j-1}=\lambda$, otherwise the operations $n_{j-1}\leftarrow-1$ and $n_{j-3}\leftarrow +1$ on $X$ produces a sequence that satisfies all conditions (B1) - (B6) and beats sequence $X$. Next we consider $n_{j-2}$. We note that $n_{j-2}\geq \lambda +1$, otherwise the operations $n_{j-1}\leftarrow-1$ and $n_{j}\leftarrow +1$ on $X$ produces a sequence that satisfies all conditions (B1) - (B6) and beats sequence $X$. Similarly, $n_{j-3}\geq \lambda +1$, otherwise the operations $n_{j-2}\leftarrow-1$ and $n_{j}\leftarrow +1$ on $X$ produces a sequence that satisfies all conditions (B1) - (B6) and beats sequence $X$. Hence, for $j\geq 5$ we have that $n_{j-3}\geq \lambda +1$, $n_{j-2}\geq \lambda +1$, $n_{j-1}= \lambda$, $n_{j}=1$, $n_{j+1}= \lambda$ and $n_{j+2}\geq \lambda -1$. If on the contrary we have $j=6$ or even $j\geq 7$, then, in addition to the values $n_{j-3}$, $n_{j-2}$, $n_{j-1}$, $n_{j}$, $n_{j+1}$, and $n_{j+2}$ we have the following.\\
{\sc case 4.2.1:} $j=6$.\\
{\sc claim 4.2.1.1:} $n_{j-4}= \lambda -1$.\\
Suppose to the contrary that $n_{j-4}\neq \lambda -1$. Clearly $n_{j-4}\geq \lambda -1$ by (B5). Hence we just have to show the contrary that $n_{j-4}\geq \lambda$. Now, if $X^*$ is a sequence obtained from $X$ by applying $n_{j-4}\leftarrow-1$ and $n_{j-2}\leftarrow+1$ clearly, $X^*$ satisfies all conditions (B1) - (B6). Then $f(X^*)=f(X)$ and $g(X^*)=g(X)+2(n_{j-4}-(\lambda -1))>g(X)$. It follows that the sequence $X^{*}$ beats sequence $X$, a contradiction. Hence $n_{j-4}=\lambda -1$.\\
{\sc case 4.2.2:} $j\geq 7$.\\
{\sc claim 4.2.2.1:} $n_{j-4}= 2$.\\
We note that $n_{j-5}\leq \lambda $, otherwise the operations $n_{j-4}\leftarrow-1$ and $n_{j-2}\leftarrow +1$ on $X$ produces a sequence that satisfies all conditions (B1) - (B6) and beats sequence $X$. Now, suppose to the contrary that $n_{j-4}\neq 2$. Clearly $n_{j-4}$ must be at least $2$ by the choice of $j$. So, we just need to suppose that $n_{j-4}\geq 3$ Now, if $X^*$ is a sequence obtained from $X$ by applying $n_{j-4}\leftarrow-1$ and $n_{j-2}\leftarrow+1$ clearly, $X^*$ satisfies all conditions (B1) - (B6). Then $f(X^*)=f(X)+(\lambda -n_{j-5})\geq f(X)$ and $g(X^*)=g(X)+2>g(X)$. It follows that the sequence $X^{*}$ beats sequence $X$, a contradiction. Hence $n_{j-4}=2$.\\
Now, fix $k<j-1$ and assume $k$ is the largest integer with $n_{k}\leq \lambda$. Let $X^*$ be the sequence obtained from $X$ by applying $n_{k+1}\leftarrow-1$ and $n_{k+3}\leftarrow+1$ clearly, $X^{*}$ satisfies all conditions (B1) - (B6). Then $f(X^*)=f(X)+(n_{k+4}-n_{k})\geq f(X)$ and $g(X^*)=g(X)+2>g(X)$. It follows that the sequence $X^{*}$ beats sequence $X$, a contradiction.\\
Hence, $n_i=1$ only if $i=0$ or $i=d$, as desired.\\[1mm]
{\sc Claim 5:} $\left(n_0,n_1,n_2\right)=\left(1,\lambda, \lambda-1\right)$.\\
We note that $n_0= 1$ is given by (B1), and $n_1= \lambda$ follows by Claim 1 for $i-1=0$. Hence, we are left to show that $n_{2}=\lambda -1$.\\
{\sc Claim 5.1:} $n_2=\lambda-1$.\\
We know that $n_{2}\geq\lambda-1$ by considering $i=1$ and condition (B5). We are thus left to show that $n_2\leq \lambda -1$. Suppose to the contrary that $n_2\geq \lambda$. To proceed we consider two cases on $n_3$. Since $d\geq 6$, we know by Claim 4 that $n_3\geq 2$.\\
{\sc Case 5.1.1:} $n_3=2$.\\
Since $d\geq 6$ and by (B4) we have that $n_4\geq 2$. Let $j\geq 3$ be the smallest integer such that $n_j\geq 3$. It follows from (B4) that for all $1\leq i\leq j-1$ we have $n_i\geq 2$. Also, by Claim 3 we have that $(n_{j+1}, n_{j+2}, n_{j+3})= (n_{j+1}, \lambda, 1)$, where $n_{j+1}\geq 3$. Now, if $X^*$ is the sequence obtained from $X$ by applying $n_2\leftarrow-(n_2+1-\lambda)$ and $n_{j+1}\leftarrow+(n_2+1-\lambda)$. Then clearly $X^*$ satisfies (B1), (B2), (B3) and (B6). We now prove that $X^*$ satisfy (B4) and (B5). It follows that since $n'_2=\lambda -1$, $n'_{j+1}>n_{j+1}$ and $n_1, n_3\geq 2$, $n'_{j}=n_{j}$, $n'_{j+2}=n_{j+2}$ then (B4) holds for all $i$ since $n_2$ and $n_{j+1}$ are the only values of $X$ that have changed. We now show that (B5) holds. Clearly (B5) holds for all $i$ except possibly $i=1,3, n_{j}, n_{j+2}$. Clearly for $i=n_{j}, n_{j+2}$ we have that (B5) holds since $n'_{j+1}> n_{j+1}$. Since $n'_{2}=\lambda -1$ any positive integer for $n_0$ or $n_4$ leaves (B5) valid. Hence (B5) holds for all $i$. And $f(X^*)=f(X)+(n_{j}-2)(n_2+1-\lambda)> f(X)$, it then follows that the sequence $X^{*}$ beats sequence $X$, a contradiction. Hence, this case cannot occur.\\
{\sc Case 5.1.2:} $n_3\geq 3$.\\
We prove separately for different choices of $n_4$. Since $d\geq 6$, we know by ($P3$) that $n_4\geq 2$.\\
{\sc Case 5.1.2.1:} $n_4=2$.\\
As in Case 5.1.1, define $j$ to be the smallest integer such that $j\geq 3$ and $n_j\geq 3$. Since $d\geq 6$ and by Claim 3 we have that $(n_{j-1}, n_{j}, n_{j+1}, n_{j+2}, n_{j+3})=(2, n_{j}, n_{j+1}, \lambda, 1)$. Now, if $X^*$ is the sequence obtained from $X$ by simultaneously applying $n_2\leftarrow-(n_2+1-\lambda)$ and $n_{j+1}\leftarrow+(n_2+1-\lambda)$, $n_3\leftarrow-(n_3-2)$ and $n_j\leftarrow+(n_3-2)$. Then clearly $X^*$ satisfies (B1), (B2), (B3) and (B6). We now prove that $X^*$ satisfy (B4) and (B5). Clearly (B4) holds for all $i$ except possibly $i=1,2,3,j-1, j,j+1$ since only the values $n_2, n_3, n_{j}$ and $n_{j+1}$ have changed. Clearly, (B4) holds for all $n_i$ where $i=j-1, j, j+1$ since $n'_{i}\geq n_{i}$ for all $i=j-1, j, j+1$. For $i=1$, since $n_1=\lambda$, any positive integer for $n_2$ leaves (B4) valid. Now, consider $i=2$, it is clear that (B4) holds since $n'_2=\lambda -1\geq 2$ and $n'_3=2$. Also, for $i=3$, (B4) holds since $n'_3, n'_{4}=2$, and hence holds for all $i$. We now show that (B5) holds. Clearly (B5) holds for all $i$ except possibly $i=1,2,3,4, j-1, j,j+1, j+2$ since only the values $n_2, n_3, n_{j}$ and $n_{j+1}$ have changed. Similarly, (B5) holds for all $n_i$ where $i=j-1, j, j+1, j+2$ since $n'_{i}\geq n_{i}$ for all $i=j-1, j, j+1, j+2, j+3$. For $i=1,2,3,4$, we have that $n'_0+n'_2=1+(\lambda -1)=\lambda$, $n'_1+n'_3=\lambda+2>\lambda$, $n'_2+n'_4=(\lambda -1)+2>\lambda$ and $n'_3+n'_5=2+n'_{5}\geq\lambda$ respectively. Hence, (B5) holds for all $i$. And, since $f(X^*)\geq f(X)$, $g(X^*)=g(X)+(n_2-\lambda)(j-1)+(n_3-1)(j-3)>g(X)$. It follows that the sequence $X^{*}$ beats sequence $X$, a contradiction. Hence, this case cannot occur.\\
{\sc Case 5.1.2.2:} $n_4\geq 3$.\\
Since $d\geq 6$ and $n_5\geq 2$ by Claim 4.\\ Since $(n_0,n_1,n_2,n_3,n_4,n_5)= (1,\lambda,n_2,n_3,n_4,n_5)$, where $n_{2}\geq \lambda$, $n_{3}, n_{4}\geq 3$ we note that if $n_5\leq \lambda -1$ then $n_6\geq 2$ by condition (B5) and hence the sequence resulting from the operation $n_4\leftarrow-1$ and $n_2\leftarrow+1$, satisfies all conditions (B1) - (B6) and beats sequence $X$. Hence, we assume $n_{5}\geq \lambda$. But, if $X^*$ is the sequence obtained from $X$ by applying $n_2\leftarrow-1$ and $n_4\leftarrow+1$, then clearly $X^*$ satisfies (B1), (B2), (B3) and (B6). We now prove that $X^*$ satisfy (B4) and (B5). It follows that since $n'_2\geq\lambda -1$, $n'_{5}>n_{5}$ and $n_1, n_3\geq 2$, $n'_{4}=n_{4}$, $n'_{6}=n_{6}$ then (B4) holds for all $i$ since $n_2$ and $n_{4}$ are the only values of $X$ that have changed. We now show that (B5) holds. Clearly (B5) holds for all $i$ except possibly $i=1,3$ and $n_{5}$. Clearly for $i=n_5$ we have that (B5) holds since $n'_{4}> n_{4}$. Since $n'_{2}\geq\lambda -1$ any positive integer for $n_0$ or $n_4$ leaves (B5) valid. Hence (B5) holds for all $i$. And $f(X^*)=f(X)+(n_5-\lambda)\geq f(X)$ and $g(X^*)=g(X)+2>g(X)$. It follows that the sequence $X^{*}$ beats sequence $X$, a contradiction. Hence, this case cannot occur.\\
Therefore, $(n_0,n_1,n_2)=(1,\lambda,\lambda -1)$, as desired.\\[1mm]
{\sc Claim 6:} If $\lambda=3,4$ and $3\leq i\leq d-4$, then $n_i=2$.\\
We know from Claim 4 that $n_i\geq 2$ for all $i=1,2,\ldots , d-1$. Now we left to show that $n_i\leq 2$ for all $i=1,2,\ldots, d-1$ to complete the proof for Claim 3.\\
Suppose to the contrary that for some $i\in\{1,2,\ldots, d-1\}$, we have $n_i\geq 3$. Let $j$ be the smallest integer such that $n_j\geq 3$ and $1\leq j\leq d-1$. Then, by Claim 3 we know that $j=d-3$, hence $n_i=2$ for all $i=3, 4,\dots, j-1$. Since $j=d-3$ we have that $i=j-1=d-3-1=d-4$. Therefore, $n_i=2$ for all $i=3, 4,\ldots, d-4$, as desired.\\[1mm]
{\sc Claim 7:} If $\lambda =2$ and $d$ odd then\\
$\left(n_{d-3}, n_{d-2}, n_{d-1}, n_{d}\right)=\left(\big\lfloor\frac{2n-3d+3}{4}\big\rfloor,\big\lceil\frac{2n-3d+3}{4}\big\rceil, 2,1\right)$.\\
Since $n_d=1$, we know that $n_{d-1}\geq \lambda$ by (B4). We first show that $n_{d-1}\leq \lambda$. 
It follows from Claim 1 and Claim 2 that for $0\leq i\leq d-4$ we have $n_i=1$ if $i$ is even, and $n_{i}=2$ if $i$ is odd.
Now, since $(n_{d-4}, n_{d-3}, n_{d-2}, n_{d-1}, n_{d})=(2, n_{d-3}, n_{d-2}, n_{d-1}, 1)$, we note that $n_{d-1}\leq \lambda$ otherwise if $n_{d-1}>\lambda$, then the sequence obtained from $X$ by applying $n_{d-1}\leftarrow-(n_{d-1}-\lambda)$ and $n_{d-3}\leftarrow+(n_{d-1}-\lambda)$ beats $X$. Thus $n_{d-1}\leq \lambda$, which implies $n_{d-1}=\lambda$. Now, we have that $(n_{d-4}, n_{d-3}, n_{d-2}, n_{d-1}, n_{d})=(2, n_{d-3}, n_{d-2}, 2, 1)$. Since $\lambda = 2$, we know by (B4) that $n_{d-3}+n_{d-2}\geq 3$. Also, $n_{d-3}+n_{d-2}=\frac{1}{2}\left(2n-3d+3\right)$.\\
{\sc Claim 7.1:} $0\leq n_{d-2}-n_{d-3}\leq 1$\\
Clearly, $n_{d-2}-n_{d-3}\geq 0$ since otherwise, if $n_{d-3}>n_{d-2}$, swapping the values of $n_{d-3}$ and $n_{d-2}$ would yield a sequence that satisfies (B1) - (B6) and has the same $f$-value, but a greater $g$-value, a contradiction. Hence it is only left to show that $n_{d-2}-n_{d-3}\leq 1$. Suppose to the contrary that $n_{d-2}-n_{d-3}\geq 2$. Let $X^*$ be the sequence obtained from $X$ by applying $n_{d-2}\leftarrow-1$ and $n_{d-3}\leftarrow+1$. Then clearly $X^*$ satisfies (B1), (B2), (B3) and (B6). We now show that (B4) holds in $X^*$. We have that (B4) holds for all $i$ except possibly when $i=d-4, d-3, d-2$ since only the values $n_{d-3}$ and $n_{d-2}$ have changed. Now, since $n_{d-4}'=n_{d-4}=2$ then any positive integer for entry $d-3$ will make (B4) hold. But $n_{d-3}'>n_{d-3}$ hence (B4) holds for $i=d-4$. For $i=d-3$, we note that $n_{d-2}\geq n_{d-3}+2\geq 3$, hence the operation leaves $n_{d-2}'\geq 2$ and $n_{d-3}'=n_{d-3}+1\geq 2$, which implies (B4) holds for $i=d-3$. Finally, for $i=d-2$, since $n_{d-1}'=n_{d-1}=\lambda$ and $n_{d-2}'= n_{d-2}-1\geq 2$ we have that (B4) holds for $i=d-2$, hence (B4) holds for all $i$. Next, we now show that (B5) holds. Since only the values $n_{d-3}$ and $n_{d-2}$ have changed, condition (B5) holds for all $i$ except possibly $i=d-4, d-3, d-2, d-1$. Now, by (B3) and since $\lambda =2$, it is clear that (B5) holds for $i=d-4, d-3, d-2, d-1$, hence (B5) holds for all $i$. Since  $f(X^*)=f(X)+n_{d-2}-(n_{d-3}+1)>f(X)$, we have that $X^*$ beats $X$, a contradiction, thus proving the claim to be true.\\[1mm]
{\sc Claim 8:} Let $\lambda =2$ and $d$ even. If
\begin{enumerate}
\item[(i)] $n-\frac{3}{2}d-5=0$ then $\left(n_{d-3}, n_{d-2}, n_{d-1}, n_{d}\right)=\left(2,1,2,1\right)$,
\item[(ii)] $n-\frac{3}{2}d-6\geq 0$ then $\left(n_{d-3}, n_{d-2}, n_{d-1}, n_{d}\right)=\left(2,\big\lceil\frac{2n-3d-4}{4}\big\rceil, \big\lfloor\frac{2n-3d-4}{4}\big\rfloor,1\right)$.
\end{enumerate}
Since $n_d=1$, we know that $n_{d-1}\geq \lambda$ by (B4). It follows from Claim 1 and Claim 2 that for $0\leq i\leq d-3$ we have $n_i=1$ if $i$ is even, and $n_i=2$ if $i$ is odd. Now, since $(n_{d-3}, n_{d-2}, n_{d-1}, n_{d})=(2, n_{d-2}, n_{d-1}, 1)$, and $\lambda =2$, we know by (B4) that $n_{d-2}+n_{d-1}\geq 3$. Also, $n_{d-2}+n_{d-1}=\frac{1}{2}\left(2n-3d-4\right)$. Now, since $n_{d-3}>n_{d}$ and $n_{d-1}\geq \lambda$, we have $n_{d-1} \geq 2$ by (B4), and $n_{d-2}\geq 1$. Hence $n_{d-1}+n_{d-2} \geq 3$. If $n_{d-2}+n_{d-1}=3$, or equivalently if $n=\frac{3}{2}d+5$, then $n_{d-2}=1$ and $n_{d-1}=2$.\\ 
If $n_{d-2}+n_{d-1} \geq 4$,  or equivalently if $n\geq \frac{3}{2}d+6$, then $0\leq n_{d-2}-n_{d-1} \leq 1$. Indeed, if $n_{d-2}-n_{d-1} <0$ then the operations $n_{d-1} \leftarrow -1$ and $n_{d-2} \leftarrow +1$ yield a sequence that satisfies (B1) - (B6) and beats $X$, and  if $n_{d-2}-n_{d-1} >1$ then the operations $n_{d-1} \leftarrow +1$ and $n_{d-2} \leftarrow -1$ again yield a sequence that satisfies (B1) - (B6) and beats $X$. We conclude that $n_{d-2}=\lceil \frac{2n-3d-4}{4} \rceil$ and $n_{d-1} =\lfloor \frac{2n-3d-4}{4} \rfloor$, as desired.\\[1mm]
{\sc Claim 9:} 
\begin{itemize} 
\item Let $\lambda =3$. If
\begin{itemize}
\item $n$ odd, then $\left(n_{d-3}, n_{d-2}, n_{d-1}, n_{d}\right)=\left(\frac{n-2d+1}{2},\frac{n-2d+3}{2}, 3,1\right)$,
\item $n$ even and $n-2d-2=0$, then $\left(n_{d-3}, n_{d-2}, n_{d-1}, n_{d}\right)=\left(2,2,3,1\right)$,
\item $n$ even and $n-2d-4\geq 0$, then

\hfill $\left(n_{d-3}, n_{d-2}, n_{d-1}, n_{d}\right)=\left(\frac{n-2d}{2},\frac{n-2d+4}{2}, 3,1\right)$.
\end{itemize}
\item Let $\lambda =4$. If
\begin{itemize}
\item $n$ odd, then $\left(n_{d-3}, n_{d-2}, n_{d-1}, n_{d}\right)=\left(\frac{n-2d-3}{2},\frac{n-2d+1}{2}, 4,1\right)$,
\item $n$ even and $n-2d-6=0$, then $\left(n_{d-3}, n_{d-2}, n_{d-1}, n_{d}\right)=\left(2,3,4,1\right)$,
\item $n$ even and $n-2d-8\geq 0$, then 

\hfill $\left(n_{d-3}, n_{d-2}, n_{d-1}, n_{d}\right)=\left(\frac{n-2d-20}{2},\frac{n-2d+2}{2}, 4,1\right)$.
\end{itemize}
\end{itemize}
We note that $n_{d}=1$ is given by (B6), hence we only left to show that $n_{d-1}=1$. Since $n_d=1$, it follows by condition (B4) that $n_{d-1}\geq \lambda$. We claim that $n_{d-1}$ is not greater than $\lambda$. Suppose to the contrary that it is. We know by Claim 4 that $n_{d-4} >1$, hence if $X^*$ is the sequence obtained from $X$ by applying $n_{d-1}\leftarrow-1$ and $n_{d-3}\leftarrow+1$. Then clearly $X^*$ satisfies (B1), (B2), (B3) and (B6). We now show that (B4). Now, clearly $i=d-4, d-3, d-2, d-1$ are the only positions in $X^*$ where possibly condition (B4) may fail since the only values $n_{d-3}$ and $n_{d-1}$ have changed. Since $n'_{d-3}>n_{d-3}$, $n'_{d-4}=n_{d-4}$ and $n'_{d-2}=n_{d-2}$, it is clear that for $i=d-4, d-3$ (B4) holds. And for $i=n_{d-2}$, we have that $n'_{d-2}n'_{d-1}=n_{d-2}(n_{d-1}-1)\geq (\lambda -1)(\lambda)>\lambda$, also, for $i=n_{d-1}$, we have that $n'_{d-1}n'_{d}=(n_{d-1}-1)n_{d}\geq \lambda$. Hence, (B4) holds for all $i\in \{0,1,\ldots, d-1\}$. We now show that (B5) holds. Since only the values $n_{d-3}$ and $n_{d-1}$ have changed, condition (B5) holds for all $i$ except possibly $i=d-4, d-2, d$. Since $n'_{d-3}\geq n_{d-3}$ and $n'_{d-5}=n_{d-5}$, it is clear that for $i=d-5$ (B5) holds. Also, since $n'_{d-3}+n'_{d-1}=n_{d-3}+n_{d-1}$, it is trivial that (B5) holds for $i=d-2$. And finally, for $i=d$ we have $n'_{d-1}+n'_{d+1}=(n_{d-1}-1)+0\geq \lambda$, hence condition (B5) holds for $i=d$ and thus holds for all $i$. Since $f(X^*)=f(X)+1>f(X)$ we have that $X^*$ beats $X$, a contradiction. Therefore $n_{d-1}=\lambda$.\\
Since $n_{d-4}=2$ by Claim 6 and $\lambda=3,4$, we have that $n_{d-3}\geq 2$ in order to satisfy condition (B4). Now, for $n_{d-2}$, it is clear that $n_{d-2}\geq \lambda -1$ by condition (B5) and (B6). Since $(n_0, n_1,n_2)=(1,\lambda, \lambda -1)$, $(n_{d-1}, n_d)=(\lambda,1)$ and $n_{i}=2$ for all $i=3,4,\ldots, d-4$, we know that $n_{d-3}+n_{d-2}=n-2d-3\lambda+11$. Since $n_{d-1}=\lambda$ is greater than $n_{d-4}=2$, we have that $n_{d-2}\geq n_{d-3}$, otherwise the operations $n_{d-3}\leftarrow-1$ and $n_{d-2}\leftarrow+1$ on $X$ produces a sequence that satisfies all conditions (B1) - (B6) and beats sequence $X$.\\
{\sc Claim 9.1:} $0\leq n_{d-2}-n_{d-3}\leq \lambda -1$.\\
We first show that $n_{d-2}-n_{d-3}=0$ if and only if $\lambda =3$ and $n-2d-3\lambda+11=4$. One direction of this claim is trivial since $n_{d-2}\geq 2$ by (B5) and $i=d-1$, also $n_{d-3}\geq 2$ by (B4), $i-d-4$ and since $n_{d-4}=2$. Next we argue that $n_{d-2}=n_{d-3}$ happens only when $\lambda=3$ and $n-2d-3\lambda+11=4$. We will consider the two cases, $\lambda\neq 3$ and $n-2d-3\lambda+11\neq 4$ and show that for each one of the cases we have that $n_{d-2}\neq n_{d-3}$. \\
{\sc CASE 9.1.1:} $\lambda\neq 3$.\\
Then $\lambda$ can only be $4$, also since $n_{d-2}\geq 3$ by (B5) and $i=d-1$, and $n_{d-3}\geq 2$ by (B4), $i=d-4$ and $n_{d-4}=2$ we have that $n-2d-3\lambda +11\geq 5$. Since $n_{d-1}=4$ and $n_{d-4}=2$ we have that $n_{d-2}$ must be greater than $n_{d-3}$, otherwise the operations $n_{d-3}\leftarrow-1$ and $n_{d-2}\leftarrow +1$ on $X$ produces a sequence that satisfies all conditions (B1) - (B5) and beats sequence $X$.\\
{\sc CASE 9.1.2:} $n-2d-3\lambda+11\neq 4$.\\
By (B5) applied on $i=d-1$ and (B4) applied on $i=d-4$ we have that $n-2d-3\lambda+11\neq 4$ implies that $n-2d-3\lambda+11\geq 5$. Since we know that $n_{d-1}\geq \lambda -1$ and since $n_{d-4}=2$ we argue in a similar way as in Case 9.1.1 that for both $\lambda =3, 4$ $n_{d-2}>n_{d-3}$, otherwise the operations $n_{d-3}\leftarrow-1$ and $n_{d-2}\leftarrow +1$ on $X$ produces a sequence that satisfies all conditions (B1) - (B6) and beats sequence $X$.\\
Now the two cases tell us that unless $\lambda=3$ and $n-2d-3\lambda+11=4$ we will always have that $n_{d-2}-n_{d-3}>0$, as desired.\\
To show that $n_{d-2}-n_{d-3}\leq \lambda -1$. Suppose to the contrary that $n_{d-2}-n_{d-3}\geq \lambda$. Now, since $n_{d-2}\geq \lambda -1$ and $n_{d-3}\geq 2$ we have that $n_{d-2}> \lambda -1 +n_{d-3}\geq \lambda +1$. Since $n_{d-4}<n_{d-1}$, let $X^*$ be the sequence obtained from $X$ by applying $n_{d-3}\leftarrow-1$ and $n_{d-2}\leftarrow+1$. Then clearly $X^*$ satisfies (B1), (B2), (B3) and (B6). Since $n_{d-3}\geq 2$ we have that $n_{d-2}\geq \lambda +1$. We now show that (B4) holds in $X^*$. Clearly $i=d-4, d-3, d-2$ are the only positions in $X^*$ where possibly condition (B4) may fail since the only values $n_{d-3}$ and $n_{d-2}$ have changed. Since $n'_{d-2}\geq \lambda >n_{d-2}$ and $n'_{d-1}=n_{d-1}$, it is clear that for $i=d-2$ (B4) holds. And for $i=n_{d-3}$, we have that $n'_{d-3}n'_{d-2}\geq (n_{d-3}+1)(n_{d-2}-1)\geq 3\lambda>\lambda$, also, for $i=n_{d-4}$, we have that $n'_{d-4}n'_{d-3}=n_{d-4}(n_{d-3}+1)\geq 6>\lambda$, for $\lambda=3,4$. Hence, (B4) holds for all $i\in \{0,1,\ldots, d-1\}$. We now show that (B5) holds. Since only the values $n_{d-3}$ and $n_{d-2}$ have changed, condition (B5) holds for all $i$ except possibly $i=d-4, d-3, d-2, d-1$. Since $n'_{d-2}\geq \lambda$, $n'_{d-4}=n_{d-4}$ and $n'_{d}=n_{d}$, it is clear that for $i=d-3$ and $i=d-1$ (B5) holds. Also, since $n'_{d-5}+n'_{d-3}=n_{d-5}+(n_{d-3}+1)\geq 5>\lambda$, we have that $i=d-4$ satisfies (B5), similarly for $i=d-2$, we have that $n'_{d-3}+n'_{d-1}=(n_{d-3}+1)+n_{d-1}\geq \lambda +3>\lambda$. Hence condition (B5) holds for all $i$. Since $f(X^*)=f(X)+(n_{d-2}-n_{d-3})-(\lambda -1)>f(X)$, we have that $X^*$ beats $X$, a contradiction.\\
This proves the claim and hence completes the proof for the lemma.
\hfill$\Box$\\

\noindent Recall that for a finite sequence $X=\left(x_0, x_1,\cdots, x_d\right)$ of positive integers we define the graph $G(X)$ by \[G(X)=\overline{K}_{x_0} +\overline{K}_{x_1}+\cdots +\overline{K}_{x_d}.\]

\begin{theorem}\label{thm4.7}
Let $G$ be a $\lambda$-edge-connected, bipartite graph of order $n$ and diameter $d$, where $d\geq 6$, $\lambda =2,3,4$ and $n \geq 10$. Then 
\[ m(G) \leq m(G(X_{n,d}^{\lambda})). \]
\end{theorem}

{\bf Proof:} 
Let $n, d$ and $\lambda$ be fixed. Assume $G$ is a $\lambda$-edge-connected bipartite graph of order $n$ and diameter $d$, and $u$ is a vertex of eccentricity $d$. Let $X(u)$ be the distance degree of $u$. It follows by Proposition \ref{prop4.3} that 
\[m(G) \leq f(X(u)).\]
Since $X(u)$ satisfies (B1) - (B6), and since by Lemma \ref{lem4.6} sequence $X(u)$ does not beat $X_{n,d}^{\lambda}$, we have 
\[ f(X) \leq f(X_{n,d}^{\lambda})\quad{\rm for}\quad d\geq 6.\]
Similarly, from Lemma \ref{lem4.5}, we have
\[ f(X) \leq f(X_{n,d}^{\lambda})\quad{\rm for}\quad d=3,4,5.\]
If $u$ is the vertex of $G(X_{n,d}^{\lambda})$ contained in $K_{x_0}$, then 
the distance degree of $u$ is $X_{n,d}^{\lambda}$, so 
\[ f(X_{n,d}^{\lambda}) \leq m(G(X_{n,d}^{\lambda})). \] 
Combining the all these (in)equalities yields the theorem.
\hfill$\Box$\\

For small diameter, the maximum size of a $\lambda$-edge-connected bipartite graph is given as in the following theorem. 

\begin{theorem}\label{lem4.5}
Let $\lambda \in \{2,3,4\}$, let $n\in \mathbb{N}$, $d\in \{3,4,5\}$, and let  $G$ be a $\lambda$-edge-connected bipartite graph of order $n$ and diameter $d$ that has maximum size among all such graphs. Then \\
(a) if $d=3$, then $m(G) \leq m\left(\overline{K}_1 + \overline{K}_{\big\lfloor\frac{n-2}{2}\big\rfloor} + \overline{K}_{\big\lceil\frac{n-2}{2}\big\rceil} + \overline{K}_1\right)$, \\
(b) if $d=4$, then 
\[m(G) \leq \left\{ \begin{array}{c}
m\left(\overline{K}_1 + \overline{K}_{\lambda} + \overline{K}_{n-2\lambda-2} + \overline{K}_{\lambda} + \overline{K}_1\right)\\
{\rm if}\quad 3\lambda +1\leq n\leq 4\lambda +1,\\
\\
m\left(\overline{K}_1 + \overline{K}_{\lambda} + \overline{K}_{\frac{n-4}{2}} + \overline{K}_{\frac{n-2\lambda}{2}} + \overline{K}_1\right)\\
{\rm if}\quad n\geq 4\lambda +2\quad {\rm and}\quad n\quad {\rm even},\\
\\
m\left(\overline{K}_1 + \overline{K}_{\lambda} + \overline{K}_{\frac{n-3}{2}} + \overline{K}_{\frac{n-\lambda -1}{2}} + \overline{K}_1\right)\\
{\rm if}\quad n\geq 4\lambda +2\quad {\rm and}\quad n\quad {\rm odd},
\end{array} \right.\]
(c) if $d=5$, then 
   \[m(G) \leq m\left(\overline{K}_1 + \overline{K}_{\lambda} + \overline{K}_{\big\lfloor\frac{n-2(\lambda +1)}{2}\big\rfloor} + \overline{K}_{\big\lceil\frac{n-2(\lambda +1)}{2}\big\rceil}  + \overline{K}_{\lambda}+ \overline{K}_1\right).\] 
\end{theorem}

{\bf Proof:} 
Let $d=3$ and $v$ be a peripheral vertex of $G$ and $X(v)=(n_0,n_1,n_2,n_3)$ a distance degree of $v$. Then we have that $G(X)$ is bipartite. We assume by Proposition \ref{prop3.3} that $G=G(X)$. Hence it should suffice to show the integers $n_0, n_1, n_2, n_3$. Following from Proposition \ref{prop4.4} we have that $\overline{K}_{n_0}=\overline{K}_1$ and $\overline{K}_{n_3}=\overline{K}_1$. Since by (B4) we have $n_1, n_2\geq \lambda$, we have $n \geq 2\lambda +2$. We also have $0 \leq n_2-n_1 \leq 1$. Indeed, if $n_2-n_1 <0$ then it is easy to verify that the size of $\overline{K}_1 + \overline{K}_{n_1-1} + \overline{K}_{n_2+1} + \overline{K}_1$ is not less than the size of $\overline{K}_1 + \overline{K}_{n_1} + \overline{K}_{n_2} + \overline{K}_1$. Similarly, if $n_2-n_1 >1$ then it is easy to verify that the size of $\overline{K}_1 + \overline{K}_{n_1+1} + \overline{K}_{n_2-1} + \overline{K}_1$ is not less than the size of $\overline{K}_1 + \overline{K}_{n_1} + \overline{K}_{n_2} + \overline{K}_1$. Hence (a) follows. 
If $d\in \{4,5\}$, then we may assume that $n_1=\lambda$ since otherwise, 
if $n_1 \geq \lambda+1$, then $\overline{K}_1 + \overline{K}_{n_1-1} + \overline{K}_{n_2+1} + \overline{K}_{n_3} + \overline{K}_{1}$ 
(or $\overline{K}_1 + \overline{K}_{n_1-1} + \overline{K}_{n_2+1} + \overline{K}_{n_3} + \overline{K}_{n_4} + \overline{K}_{1}$) 
has size not less than $\overline{K}_1 + \overline{K}_{n_1} + \overline{K}_{n_2} + \overline{K}_{n_3} + \overline{K}_{1}$ 
(or $\overline{K}_1 + \overline{K}_{n_1} + \overline{K}_{n_2} + \overline{K}_{n_3} + \overline{K}_{n_4} + \overline{K}_{1}$). 
Similarly for $d=5$ or $d=4$ and $3\lambda +1\leq n\leq 4\lambda +1$ we may assume that $n_{d-1}=\lambda$. If $d=4$ and $3\lambda +1\leq n\leq 4\lambda +1$ we 
get that $m(G) \leq m(K_1 + K_{\lambda} + K_{n-2\lambda-2} + K_{\lambda} + K_1)$. Otherwise, let $d=4$ and $n\geq 4\lambda +2$. We also know that by (B4) and (B5) we have $n_2\geq \lambda -1$ and $n_3\geq \lambda$. We also have $\lambda -2 \leq n_2-n_3 \leq \lambda -1$. Indeed, if $n_2-n_3 <\lambda -3$ then it is easy to verify that the size of $\overline{K}_1 + \overline{K}_{\lambda} +\overline{K}_{n_2+1} + \overline{K}_{n_3-1} + \overline{K}_1$ is not less than the size of $\overline{K}_1 +\overline{K}_{\lambda}+ \overline{K}_{n_2} + \overline{K}_{n_3} +\overline{K}_1$. Similarly, if $n_2-n_3 >\lambda -1$ then it is easy to verify that the size of $\overline{K}_1 + \overline{K}_{\lambda} + \overline{K}_{n_2-1} + \overline{K}_{n_3+1} +\overline{K}_1$ is not less than the size of $\overline{K}_1 + \overline{K}_{\lambda} +\overline{K}_{n_2} + \overline{K}_{n_3} + \overline{K}_1$. Hence (b) follows, Let $d=5$. Since by (B5) we have $n_2, n_3\lambda -1$, we have $n \geq 4\lambda$ and in addition, if $\lambda=2$ we have by (B4) that at least one of either $n_2$ or $n_3$ is at least $2$. Hence for $\lambda=2$ we that $n> 4\lambda$. We also have $0 \leq n_3-n_2 \leq 1$. Indeed, if $n_3-n_2 <0$ then it is easy to verify that the size of $\overline{K}_1 + \overline{K}_{\lambda} \overline{K}_{n_2-1} + \overline{K}_{n_3+1} +\overline{K}_{\lambda}+ \overline{K}_1$ is not less than the size of $\overline{K}_1 +\overline{K}_{\lambda}+ \overline{K}_{n_2} + \overline{K}_{n_3} + \overline{K}_{\lambda}+\overline{K}_1$. Similarly, if $n_3-n_2 >1$ then it is easy to verify that the size of $\overline{K}_1 + \overline{K}_{\lambda} + \overline{K}_{n_2+1} + \overline{K}_{n_3-1} + \overline{K}_{\lambda} +\overline{K}_1$ is not less than the size of $\overline{K}_1 + \overline{K}_{\lambda} +\overline{K}_{n_2} + \overline{K}_{n_3} + \overline{K}_{\lambda} + \overline{K}_1$. Hence (c) follows, and the proof of the theorem is complete.
\hfill$\Box$\\

\noindent Since $G(X_{n,d}^{\lambda})$ is a bipartite graph of order $n$, diameter $d$ and edge-connectivity $\lambda$, the bounds obtained in Theorem \ref{thm4.7} and Theorem \ref{lem4.5} are sharp. Evaluating the size of $G(X_{n,d}^{\lambda})$ yields the following corollary.\\

\begin{corollary}\label{cor4.8}
Let $n, d, \lambda\in \mathbb{N}$, with $\lambda = 2,3,4$ and $d\geq 3$. If $G$ is a $\lambda$-edge-connected bipartite graph of order $n$ and diameter $d$, then
 \begin{enumerate}
\item[($a$)] If $d=3$ then
 \[ m(G)\leq  \left\{ \begin{array}{cc}
  \frac{n^2-4}{4}
        & \textrm{if $\lambda=2,3,4$ and $n\geq 2(\lambda +1)$ and $n$ is even,}\\
\frac{n^2-5}{4}
        & \textrm{if $\lambda=2,3,4$ and $n\geq 2(\lambda +1)$ and $n$ is odd.} \\           
      \end{array} \right.\] 
 
\item[($b$)] If $d=4$ then
 \[ m(G)\leq  \left\{ \begin{array}{cc}
 2(n-2\lambda -1)\lambda 
        & \textrm{if $\lambda=2,3,4$ and $3\lambda +1\leq n\leq \frac{7}{2}\lambda +2$,} \\  
  \frac{n\left(n-2\right)}{4} 
        & \textrm{if $\lambda=2,3,4$ and $n> \frac{7}{2}\lambda +2$ and $n$ is even,}\\
\frac{(n-1)^2}{4}
        & \textrm{if $\lambda=2,3,4$ and $n> \frac{7}{2}\lambda +2$ and $n$ is odd.} \\           
      \end{array} \right.\] 
      
\item[($c$)] If $d=5$ then
  \[ m(G)\leq  \left\{ \begin{array}{cc}
\frac{(n+2)^2}{4} 
        & \textrm{if $\lambda=2$, $n\geq 9$ and $n$ is even,} \\  
  \frac{(n+7)(n-3)}{4}+6
        & \textrm{if $\lambda=2$, $n\geq 9$ and $n$ is odd,}\\
\frac{(n+2)^2}{4}-3 
        & \textrm{if $\lambda=3$ and  $n\geq 12$ and $n$ is even,} \\
  \frac{(n+7)(n-3)}{4}+3
        & \textrm{if $\lambda=3$ and  $n\geq 12$ and $n$ is odd,} \\ 
\frac{(n+2)^2}{4}-8
        & \textrm{if $\lambda=4$ and  $n\geq 16$ and $n$ is even,} \\ 
  \frac{(n+7)(n-3)}{4}-2
        & \textrm{if $\lambda=4$ and  $n\geq 16$ and $n$ is odd.} \\            
      \end{array} \right. \] 
        
\item[($d$)] If $d\geq 6$ then
 \[ m(G)\leq  \left\{ \begin{array}{c}
2\left(d-3\right)+\lfloor\frac{2n-3d-4}{4}\rfloor +2\lceil\frac{2n-3d-4}{4}\rceil+\lfloor\frac{2n-3d-4}{4}\rfloor\lceil\frac{2n-3d-4}{4}\rceil \\
         \textrm{if $\lambda=2$, $n\geq 10$ and $n$ is even,} \\ 
         \\ 
  2\left(d-3\right)+2\left(\lfloor\frac{2n-3d+3}{4}\rfloor +\lceil\frac{2n-3d+3}{4}\rceil\right)+\lfloor\frac{2n-3d+3}{4}\rfloor\lceil\frac{2n-3d+3}{4}\rceil\\
         \textrm{if $\lambda=2$, $n\geq 11$ and $n$ is odd,}\\
         \\
 2(2d+1)\\
         \textrm{if $n=2d+2$ and $\lambda=3$,} \\
         \\        
  \frac{5n-13-2d}{2}+\frac{(n-2d+1)(n-2d+3)}{4}\\
         \textrm{if $n\geq 2d+3$, $n$ odd and $\lambda=3$,} \\
         \\        
  \frac{5n-12-2d}{2}+\frac{(n-2d)(n-2d+4)}{4}\\
         \textrm{if $n\geq 2d+4$, $n$ even and $\lambda=3$,} \\
         \\          
  2(2d+9)\\
         \textrm{if $n=2d+6$ and $\lambda=4$,} \\
         \\        
3(n-1)-2d+\frac{(n-2d-3)(n-2d+1)}{4}\\
         \textrm{if $n\geq 2d+7$, $n$ odd and $\lambda=4$,} \\
         \\        
 3(n-6)-2d+\frac{(n-2d-20)(n-2d+2)}{4}\\
         \textrm{if $n\geq 2d+8$, $n$ even and $\lambda=4$,} \\
         \\               
      \end{array} \right. \] 
 \end{enumerate}
and this bound is sharp.
\end{corollary}

\section*{Acknowledgment}
The author sincerely thanks Professor Peter Dankelmann for suggesting the problem, and for his valuable guidance.



\begin{thebibliography}{9}
\bibitem{al} Ali, P.; Mazorodze, J.P.; Mukwembi, S.; Vetr\'ik, T.; \textit{On size, order, diameter and edge-connectivity of graphs.} Acta Mathematica Hungarica 152.1 (2017), 11-24.

\bibitem{bo} Bollob\'{a}s, B.; {\it Graphs with given diameter and minimal degree.} Ars Combin. {\bf 2} (1976), 3-9. 

\bibitem{bol} Bollob\'{a}s, B.; {\it On Graphs with equal edge connectivity and minimum degree.} Discrete Mathematics {\bf 28} (1979), 321-323. 
 
\bibitem{c} Caccetta, L.; Smyth, W.F.; \textit{Graphs of maximum diameter.} Discrete Math.\ {\bf 102} No.\ 2 (1992), 121-141. 

\bibitem{ca} Caccetta, L.; Smyth, W.F.; \textit{Properties of edge-maximal K-edge-connected D-critical graphs.} J. Combin. Math. Combin. Comput. {\bf 2} (1987), 111-131.

\bibitem{ca2} Caccetta, L.; Smyth, W.F.; \textit{Redistribution of vertices for maximum edge count in K-edge-connected D-critical graphs.} Ars Combin. {\bf 26} (1988), 115-132. 
 
\bibitem{da} Dankelmann, P;  \textit{Size of graphs and digraphs with given diameter and connectivity constraints.} Acta Math. Hungar. 164 (2021), 178-199. 

\bibitem{dan} Dankelmann, P.; Volkmann, L.; \textit{Minimum size of a graph or digraph of given radius.} Inform. Process. Lett., {\bf 109} (2009), 971-973.    
     
\bibitem{dank} Dankelmann, P.; Domke, G. S.; Goddard, W.; Globler, P.; Hattingh, J. H.; Swart, H. C.; \textit{Maximum sizes of graphs with given domination parameters.} Discrete Mathematics {\bf 281} (2004), 137-148.  

\bibitem{danke} Dankelmann, P.; Mafunda, S.; Mallu, S.; \textit{Remoteness of graphs with given size and connectivity constraints.} Discrete Mathematics {\bf 348} (2025), 114451.

\bibitem{er} Erd\H{o}s, P.; R\'{e}nyi, A.;\textit{On a problem in the theory of graphs.} Publ. Math. Inst. Hungar. Acad. Sci. Series B. {\bf 7} (1963), 623-639.

\bibitem{erd} Erd\H{o}s, P.; R\'{e}nyi, A.; S\'{o}s, V. T.; \textit{On a problem of graph theory.}  Studia Sci. Math. Hungar {\bf 1} (1966), 215-235.

\bibitem{fu} Z. F{\"u}redi, Z.; \textit{The maximum number of edges in a minimal graph of diameter 2.} J. Graph Theory. {\bf16} (1992), 81-98.

\bibitem{Maf2020} Mafunda, S.; \textit{Aspects of distance measures in graphs and digraphs.} Ph.D. Thesis, University of Johannesburg, (2020).
   
\bibitem{mu} Mukwembi, S.; \textit{On size, order, diameter and minimum degree.} Indian J. Pure Appl. Math. {\bf 44} No. 4 (2013), 467-472.

\bibitem{or} Ore, O.; \textit{Diameters in graphs.} J. Combin. Theory {\bf 5} (1968), 75-81.

\bibitem{vi} Vizing, V.; \textit{The number of edges in a graph of given radius.} Soviet. Math. Dokl. {\bf 8} (1967), 535-536.
\end{thebibliography}
\end{document}